
\documentstyle{amsppt}
\magnification=1200
\define\Z{{\Bbb Z}}
\define\Q{{\Bbb Q}}

\define\CC{{\Cal C}}
\define\JJ{{\Cal J}}
\define\FF{{\Cal F}}
\define\A {{\Cal A}}
\define\M{{\Cal M}}
\define\C{{\Bbb C}}
\define\Pone {{\Bbb P}^1}
\define\la{\langle}
\define\ra{\rangle}
\define\refer#1{(#1)}
\define\supp{\operatorname{supp}}
\define\Aut{\operatorname{Aut}}
\define\inw{\operatorname{int}}
\newcount\headnumber
\newcount\labelnumber
\define\section{\global\advance\headnumber
by1\global\labelnumber=0{{\the\headnumber}.\ }}
\define\label{(\global\advance\labelnumber by1 \the\headnumber
.\the\labelnumber )\enspace}
\NoBlackBoxes

\topmatter
\title
On the tautological ring of $\M _g$
\endtitle
\rightheadtext{Tautological ring}
\author
Eduard Looijenga
\endauthor

\address Faculteit Wiskunde en Informatica,
Universiteit Utrecht,
P.O. Box 80.010, 3508 TA Utrecht,
The Netherlands\endaddress
\email looijeng\@math.ruu.nl\endemail

\abstract
We prove that any product of tautological classes of $\M _g$  of degree $d$
(in the Chow ring of $\M _g$ with rational coefficients)  vanishes for $d>g-2$
and is proportional to the class of the hyperelliptic locus in degree $g-2$.
\endabstract

\endtopmatter
\document

\head
\section Results
\endhead

Fix an integer $g\ge 2$ and denote by $\CC _g^n$ the moduli space
of tuples $(C,x_1,\dots ,x_n)$, where $C$ is a smooth connected projective
curve
of genus $g$ and $x_1,\dots ,x_n$ are (not necessarily distinct) points of $C$;
we also write $\M _g$ for $\CC _g^0$.
Forgetting the $i$th point defines a morphism  $\CC _g^n\to\CC _g^{n-1}$
whose  relatively dualizing sheaf is denoted by $\omega _i$ ($i=1,\dots ,n)$.
We
write   $K_i$ for the first Chern class of $\omega _i$, considered as an
element
of the Chow group $A^1(\CC _g^n)$ (with rational coefficients); for $n=1$ we
also write $K$. Our main result is:

\proclaim{\label Theorem}
Any monomial of degree $d$ in the classes $K_1,\dots ,K_n$ is a linear
combination of the classes of the irreducible components of the locus
parametrizing tuples $(C,x_1,\dots ,x_n)$ admitting a morphism $C\to\Pone$ of
degree $\le g+n$ such that the fiber over $0$ (resp.\ $\infty$) has at most
$g+n-d-1$ points (resp.\  is a singleton) and $\{ x_1,\dots ,x_n\}$ is
contained in one of these two fibers. (Hence such a class is zero when
$d>g+n-2$.) All monomials of degree $g+n-2$ are proportional to the class of
the locus parametrizing tuples $(C,x_1,\dots ,x_n)$ with $C$ hyperelliptic and
$x_1=\cdots =x_n$ a Weierstra\ss\ point.
\endproclaim

The direct image of $K^{d+1}$ in $A^d(\M _g)$ is the Mumford--Morita--Miller
{\it tautological class} $\kappa _d$. Mumford showed in his fundamental paper
\cite{4} that the subring of $A^{\bullet}(\M _g)$ generated by these classes
(the {\it tautological ring} of $\M _g$) is already generated by $\kappa
_1,\dots ,\kappa _{g-2}$. On the basis of many calculations Carel Faber has
made
the intriguing conjecture that this ring has the formal properties of the
even-dimensional cohomology ring of a projective manifold of dimension $g-2$,
i.e., satisfies Poincar\'e duality and a Lefschetz decomposition.
We offer the following support for this conjecture:

\proclaim{\label Theorem}
Any product of tautological classes that has degree $d$ is a linear combination
of the classes of the irreducible components of the locus parametrizing curves
$C$ admitting a morphism $C\to\Pone$ of degree $\le 2g-2$ totally ramified over
$\infty$ and with at most $g-1-d$ points over $0$ (hence is zero when
$d>g-2$).
All such classes of degree $g-2$ are proportional to the class of the
hyperelliptic locus.
\endproclaim

A finer analysis of our proof may well yield that $\kappa _1^{g-2}$ is a
nonzero multiple of the hyperelliptic class, but it is not known whether the
latter is actually nonzero.

The proof of the theorems uses the flag of subvarieties of $\M _g$ introduced
by
Arbarello \cite{1}, a variant of which was exploited by Diaz \cite{2} to prove
that $\M _g$ has no complete subvarieties of dimension $>g-2$. Our simple key
result
\refer{2.4} serves as a substitute for Diaz's lemma $2$ in \cite{2} and can be
used in that paper to eliminate the use of compactifications of Hurwitz schemes
(see \refer{2.8}).
The proof of the second assertion of each theorem involves an application of
the Fourier transform for abelian varieties, due to Mukai and Beauville.

\smallskip
In this paper we only consider Chow groups with respect to rational
equivalence, tensorized with $\Q$, and graded by codimension, notation:
$A^{\bullet}$. If $X$ is a variety that is smooth, or more generally, that
admits a smooth Galois covering, then there is an intersection product
$A^k(X)\otimes A^l(X)\to A^{k+l}(X)$.

\smallskip
I thank Johan de Jong for drawing my attention to the paper by Deninger--Murre
\cite{2} and for comments on a first draft.

\head
\section Proofs
\endhead

\label Let $C$ be a smooth projective curve of genus $g$ and let $D_0$ and
$D_{\infty}$  be positive divisors on $C$ that are linearly
equivalent, but whose supports are disjoint. Then there is a
finite morphism  $\pi :C\to\Pone$ such that $\pi ^*(i)=D_i$ ($i=0,\infty $). If
$p\in C$ occurs in $D_i$ with multiplicity $m_p>0$, then $\pi$
determines an isomorphism of $\C\cong T^*_i\Pone$ onto $T^*_pC ^{\otimes
m_p}$.
However, $\pi$ is not unique for it is defined up to natural action of $\C
^{\times}$ on $\Pone$. That ambiguity can be eliminated as follows.

Let $R$ denote the part of the ramification divisor of $\pi$ that lies over
$\Pone -\{ 0,\infty \}$. If $c$ denotes the number of points of
$\supp (D_0+D_{\infty})$, then the Riemann-Hurwitz formula implies that the
degree $r$ of $R$ is equal to $2g-2+c$. If $\pi _*(R)= \sum _i n_i(z_i)$, then
$\pi$ can be
normalized in such a way that  $\prod _i z_i ^{n_i}=1$. This normalization is
unique up to multiplication by an $r$th root of unity.
So for $p$ and $m_p$ as above, and $\pi$ normalized, the corresponding
generator of $T^*_pC ^{\otimes m_p}$ raised to the $r$th power
gives a {\it canonical} generator of $T^*_pC ^{\otimes m_pr}$.

This argument works just as well in families and so we obtain:

\proclaim{\label Proposition}
Let $f:\CC \to S$ be a projective family of smooth genus
$g$ curves with reduced base. Let $D_0$ and $D_{\infty}$ be positive relative
divisors on $\CC$ whose supports are disjoint and are \'etale over $S$. Suppose
that their difference  is
linearly equivalent to the pull-back of a divisor on $S$. Then for every
section
$x:S\to \CC$ of $f$ with image in the support of  $D_0+D_{\infty}$,  a suitable
positive tensor power of $x^*\omega _{\CC /S}$ is trivial.
\endproclaim

We shall use the following simple fact:

\proclaim{\label Lemma}
Let $L_1,\dots ,L_d$ be line bundles on a variety $V$ and let $V=V^0\supset
V^1\supset\cdots\supset V^d$ be a chain of closed
subvarieties such that $L_k$ is trivial on $V^{k-1}-V^k$. Then $c_1(L_1)\cdots
c_1(L_d)$ has support in $V^d$.
\endproclaim

The key result we need is:

\proclaim{\label Lemma}
Let $d$ be a positive integer and let $\{(C_t,x_t,P_t)\} _{t\in\Delta}$ be
an analytic family of triples consisting of a smooth connected projective curve
$C_t$, a point $x_t\in C_t$, and a pencil $P_t$ on $C_t$ containing $d(x_t)$.
Assume that for $t\not=0$, $P_t$ has no base points and let $R_t$ be the part
of
the ramification divisor on $C_t-x_t$ of the associated morphism $C_t\to P_t$.
If $R_0$ is the limit of $R_t$ for $t\to 0$, then the multiplicity of $x_0$ in
$R_0$ is also the multiplicity of $x_0$ as a fixed point of $P_0$.
\endproclaim \demo{Proof} Represent the family by a smooth analytic morphism
$t:\CC\to\Delta$ with section $x:\Delta \to\CC$. Extend $t$ to a chart $(z,t)$
at $x_0$ such that $z=0$ is the image of $x$ at $x_0$. In terms of these
coordinates generators of $P_t$ can be represented by $z^d$ and a holomorphic
function  $A(z,t)=\sum _{i\not= d} a_i(t)z^i$ which is divisible neither by $t$
nor by $z$.  The first index $k$ for which $a_k(0)\not= 0$ is $<d$ and is equal
to the multiplicity of $x_0$ as fixed point of $P_0$. In the domain of the
chart, the divisor $R_t$ is given by locus where the $z$-derivatives of $A$ and
$z^d$ are proportional, i.e.,  by the divisor of $\sum _{i \not= d}
(i-d)a_i(t)z^i$ ($t\not= 0$). This expression is not divisible by $z$ or $t$ so
that $R_0$ is given by  $\sum _{i\not= d} (i-d)a_i(0)z^i$. So $x_0$ occurs with
multiplicity $k$ in $R_0$. \enddemo

An immediate consequence is an amplification of a result due to Arbarello
\cite{1} and Diaz \cite{2}:

\proclaim{\label Corollary} Suppose that in the situation of \refer{2.4} there
exists an analytic section $\{ D_t\in P_t\} _{t\in\Delta}$ such that for
$t\not=
0$, $\supp (D_t)$ is disjoint with $x_t$ and has $d-r$ points, whereas
$D_0=d(x_0)$. Then $P_0$ can be written as $r(x_0) +P'$. \endproclaim

\medskip
\label If $d$ is positive integer, then we have moduli space $P(d)$ of triples
$(C,x,P)$ with $C$ a smooth projective curve of genus $g$, $x\in C$ and $P$ a
pencil on $C$ containing $(d)x$. The existence of this is clear if $d>2g-2$,
for then this is just a bundle of projective spaces of dimension $d-g-1$ over
$\CC _g$; the remaining cases $d\le 2g-2$ follows from this by simply viewing
$P(d)$ as the locus in $P(2g-1)$ parametrizing triples $(C,x,P)$ for which $x$
is a fixed point in $P$ of multiplicity $2g-1-d$. This implies that we also
have defined a moduli space $Z$ of tuples $(C,x_1,\dots ,x_n,x,D,P)$
with $C$ a smooth projective curve of genus $g$, $x_1,\dots ,x_n,x\in C$, $P$ a
pencil on $C$ containing $(n+g)x$, $D$ a degenerate member of $P$ and $\{
x_1,\cdots ,x_n\}\subset \supp (D)$.
Notice that $D$ and $x$ determine $P$ unless $D=(n+g)(x)$. The forgetful
morphism $f:Z\to\CC ^n_g$ is clearly proper.

The tuples for which $\supp (D)$ has at most $g+n-1-k$ points outside $x$
define a closed subvariety $Z^k$ of $Z$. It is clear that $Z^{n+g-1}$ can be
identified with the set of tuples $(C,x,\dots ,x,x,(n+g)x,P)$ with $P$ a pencil
through $(n+g)(x)$.

\proclaim{\label Lemma}
For $k<g+n-1$, $Z^k-Z^{k+1}$ is Zariski-open in an affine variety of pure
dimension $3g-3+n-k$ and $f^*K_i|Z^k-Z^{k+1}=0$ ($i=1,\dots ,n$).
\endproclaim
\demo{Proof}
Let $k<g+n-1$ and let $W$ be a connected component of $Z^k-Z^{k+1}$.
If $(C,x_1,\dots ,x_n,x,D,P)$ represents an element of $W$, then write
$D=m(x)+D'$ with $x\notin\supp (D')$ so that $\supp (D')$ has exactly $n+g-k$
points.  There is a finite morphism $\pi :C\to\Pone$ with $\pi ^*(0)=D'$ and
$\pi ^*(\infty )=(g+n-m)(x)$. The part of the ramification divisor $R$ of $\pi$
over $\Pone -\{ 0,\infty \}$ has by Riemann-Hurwitz degree
$2g-2+(g+n-k)=3g-2+n-k$.

The multiplicity $m$, the multiplicity of $x_i$ in $D$, and the stratum of the
diagonal stratification of $\CC ^{n+1}_g$ containing $(C,x_1,\dots ,x_n,x)$
only
depend on $W$.  So assigning to $(C,x_1,\dots ,x_n,x,P)$ the $\C
^{\times}$-orbit of $\pi _*R$ defines a flat, quasi-finite morphism from $W$ to
the quotient of a $(3g-3+n-k)$-dimensional torus by an action of the symmetric
group. So $W$ is pure of dimension $3g-3+n-k$. Proposition \refer{2.2} implies
that $f^*K_i|W$ is trivial. \enddemo

\demo{Proof of the first clause of \refer{1.1}} Let $X^k$ be the union of
irreducible components of $Z^k$ that are distinct from $Z^{n+g-1}$. (It can be
shown that $Z^{n+g-1}$ is actually an irreducible component of $Z$
and so $X^0\not= Z$.) The restriction  $f: X^0\to\CC ^n_g$ is clearly proper.
It is also surjective, because for given $(C,x_1,\dots ,x_n)$, the morphism
$$
(y,y_1,\dots ,y_{g-1})\in C^g\mapsto [-(n+g)y + 2(x_1)+\sum _{i=2}^n(x_i)+\sum
_{j=1}^{g-1}(y_j)]\in J(C)
$$
is onto. Observe that $X^{n+g-1}=\emptyset$.

We claim that $f(X^k\cap Z^{n+g-1})\subset f(X^{k+1})$. For if
$(C,x,\dots ,x,x,(n+g)x,P)$ represents an element of $X^k\cap Z^{n+g-1}$,
then by \refer{2.5}, $P$ will be of the form $(k+1)x +P'$ with $P'$ a pencil of
degree $n+g-k-1$. So $P$ has a member $\not= (n+g)(x)$ with at most $n+g-k-2$
points.

It follows that the pre-image $U^k$ of
$f(X^k)-f(X^{k+1})$ in $X^k$ is contained in $Z^k-Z^{k+1}$. In particular,
$f^*K_i|U^k=0$ for $i=1,\dots ,n$. Since $f:U ^k\to f(X^k)-f(X^{k+1})$
is proper and onto, we also have $K_i|f(X^k)-f(X^{k+1})$ =0. So by \refer{2.3},
a monomial of degree $k$ in $K_1,\dots ,K_n$ is a linear combination of
irreducible components of $f(X^k)$ of codimension $k$. One easily checks that
these components are as described in the theorem.
\enddemo

\label Since $f(X_k)-f(X_{k+1})$ admits a finite covering that is Zariski-open
in an affine variety, it cannot contain a complete curve. From this we recover
Diaz's theorem which asserts that $\CC _g^n$ does not contain a complete
subvariety of dimension $>g+n-2$.

In order to complete the proof of \refer{1.1} we need two more results, one
algebraic, one topological.

\proclaim{\label Lemma}
Let $f:\A \to S$ be a family of abelian
varieties of dimension $g$ and let $d$ be a positive integer.
Then the class of the locus $\A \la d\ra $ of points of order $d$ is
a positive multiple of the class of the zero section in $A^g(\A )$. (The
coefficient is the number of elements in $(\Z /d)^{2g}$ of order $d$.)
\endproclaim
\demo{Proof} We use the Fourier transform for abelian varieties introduced by
Mukai, developed by Beauville and extended to abelian schemes by
Deninger--Murre \cite{2}. Mukai's transform gives an (in
general inhomogeneous) isomorphism $\FF : A(\A )\to A(\hat\A )$, where
$\hat\A\to S $ is the dual family. We shall compare the images of
the two classes in $A(\hat\A )$ under $\FF$.

Let $k$ be an positive integer relative prime to $d$.
Multiplication by $k$ in $A$ maps $\A \la d\ra $ isomorphically onto itself. So
the class of $\A \la d\ra $ in $A^g(\A )$ is fixed under $k_*$.
Lemma \refer{2.18} of \cite{2} implies that then $\FF ([\A \la d\ra ])\in
A^0(\hat\A )$. Since the projection induces an isomorphism
$A^0(S )\to A^0(\hat\A )$, the lemma follows.
\enddemo

\proclaim{\label Lemma} Let $\pi :C\to\Pone$ be a covering of degree $d$
by a smooth connected curve that is totally ramified over $0$ and
$\infty$ such that the part $D$ of the discriminant in $\Pone -\{ 0,\infty\}$
is
reduced. Then there exists a disk neighborhood $B$ of
$\supp (D)$ in $\Pone -\{0,\infty\}$ such that for $p\in\partial B$, the
monodromy group
of $\pi$ over $(B-\supp (D),p)$ is a single transposition $(a',a'')$. Moreover,
if $\sigma$ is the monodromy of a simple loop in $\Pone -\inw (B)$ around $0$
based at $p$, then
$a''=\sigma ^r(a')$ for some divisor $r$ of $d$ and $\pi$ factorizes through
the covering $z\in\Pone\to z^r\in\Pone$.
\endproclaim
\demo{Proof}
We choose a base point $p\in\Pone$ outside the discriminant and we put $F:=\pi
^{-1}(p)$. By a {\it simple arc} we shall mean an embedded interval connecting
$p$ with a point of the discriminant that does not meet the discriminant along
the way. A simple arc $\alpha$
determines up to isotopy (relative $p$ and the discriminant) a simple loop
based at $p$ around a point of the discriminant and hence a
monodromy transformation $\tau _{\alpha}\in\Aut (F)$. A collection of simple
arcs that do not meet outside $p$ shall be called an {\it arc system}. Notice
that the directions of departure of the members of such a collection determine
a cyclic ordering (our preference is clockwise) of these.

We begin by fixing a simple arc $\omega$ connecting $p$ with $0$. We write
$\sigma$ for $\tau _{\omega}$; this is a $d$-cycle in $\Aut (F)$.
Any transposition $\tau$ in $\Aut (F)$ can be written $(a,\sigma ^l(a))$ for
some $l\in\{ 0,1,\dots {1\over 2}d\}$; this means that $\sigma\tau$ is the
product of two disjoint cycles of length $l$ and $d-l$. Let us call $l$ the
{\it mesh} of $\tau$.

Let $\alpha _1$ be an simple arc to a point of $\supp (D)$ that forms with
$\omega$ an arc
system and is such that $\tau :=\tau _{\alpha _1}$ has minimal mesh $r$. Write
$\sigma\tau =\sigma '\sigma ''$ with
$\sigma '$ and $\sigma ''$ disjoint cycles of length $r$ resp.\ $d-r$ and
denote by $F'$ and $F''$ the corresponding parts of $F$. Notice that $\tau
_{\alpha _1}$ interchanges some $a'\in F'$ with some $a''\in F''$.

Let $\beta$ be another simple arc to a point of $\supp (D)$ such that  $(\omega
,\alpha _1
,\beta )$ is a clockwise oriented arc system. Then $\tau _{\beta}$ cannot
commute with $\sigma ''$: if it did, then it would interchange two points of
$F'$ and would therefore have a mesh $<r$.
It may happen that $\tau _{\beta}$ commutes with $\sigma '$. But not every
choice for $\beta$ can be like this, for then $\sigma '$ would commute with the
monodromy around $\infty$ and this is impossible as the latter is a $d$-cycle.

So for some $\beta$, $\tau _{\beta}$
interchanges some $b'\in F'$ with some $b''\in F''$. If we modify $\beta$
by letting it first wind $k$ times
around the union of $\omega$ and $\alpha _1$, then its monodromy gets
conjugated by $(\sigma '\sigma '')^{\pm k}$. In this way we can arrange that
$b''=a''$. If $b'\not= a'$, then a straightforward verification shows that
$\tau _{\beta}$ would have a smaller mesh than $r$. So $b'=a'$ and hence $\tau
_{\beta}=\tau$. This argument proves more: the fact that for every integer $k$
the mesh of the $(\sigma '\sigma '')^k$-conjugate of $\tau _{\beta}$ is $\ge r$
implies that $r$ divides $d$. We put $\alpha _2:=\beta$.

We now prove with induction on $l$ that for $l\le \deg (D)$ there is an arc
system $(\alpha _1,\alpha _2,\dots ,\alpha _l)$ in clockwise cyclic order such
that $\tau _{\alpha _i}=\tau$ for $i=1,\dots ,l$.
The lemma then follows: we already showed that $r$ divides $d$, and it is easy
to see that the asserted factorization exists. So suppose we found such an arc
system $(\alpha _1,\alpha _2,\dots ,\alpha _l)$ for some $l\ge 2$.

First assume $l$ even. Then the monodromy around the union of these arcs is
equal to $\sigma$ and so the above argument yields simple arcs $\beta _1,\beta
_2$ such that $\tau _{\beta _1}=\tau _{\beta _2}$ and $(\omega ,\alpha _1,\dots
,\alpha _l,\beta _1,\beta _2)$ is an arc system in clockwise order.
Since $\tau _{\beta _i}$ does not commute with $\sigma ''$, we can modify
$\beta _1$ and $\beta _2$ by letting both go round the union of $(\omega
,\alpha _1,\dots ,\alpha _l)$ the same number of times first, to ensure that
$\tau _{\beta _1}=\tau _{\beta _2}$ moves $a''$.

If $\tau _{\beta _i}$ does not commute with $\sigma '$, then the argument above
shows that in fact $\tau _{\beta _i}=\tau$ and so we managed to take two
induction steps.

If $\tau _{\beta _i}$ does commute with $\sigma '$, then let $\beta '_i$ be
obtained from $\beta _i$ by going round $\alpha _l$ first. Then $(\omega
,\alpha _1,\dots ,\alpha _{l-1},\beta '_1,\beta ' _2,\alpha _l)$ is in
clockwise order and
$\tau _{\beta '_1}=\tau _{\beta '_2}$ interchanges an element of $F'$ with an
element of $F''$. Next modify the $\beta '_1$ and $\beta '_2$ by letting them
first encircle $(\omega ,\alpha _1,\dots ,\alpha _{l-1})$ the same number of
times as to arrange that $\tau _{\beta '_i}$ moves $a''$ (this might cause them
to meet $\alpha _l$ in a point $\not= p$). Then $\tau _{\beta '_i}=\tau$ and
hence we have made the induction step.

It remains to do the induction step for $l$ odd. That is handled in the same
way as the case $l=1$.
\enddemo

\demo{Proof of the second clause of \refer{1.1}}
Notice that $X^{n+g-2}$ parametrizes the triples $(C,x,y)$, where $C$
is smooth of genus $g$, $x,y\in C$ are distinct and $d(x)\equiv d(y)$ for some
$d\in \{ 2,\dots ,n+g\}$. By our previous discussion this defines a closed
subvariety $Y$ of $\CC ^2_g$ of pure codimension $g$. The assertion that is to
be proved will follow if we show that the classes of the irreducible components
of $Y$ are proportional in $A^g(\CC ^2_g)$. Our first business is therefore to
describe these irreducible components.

For $d\ge 2$, let $Y_d\subset \CC ^2_g$ be the locus parametrizing triples
$(C,x,y)$ for which $(x)-(y)$ has order $d$ in $J(C)$. For such
$(C,x,y)$ we have a morphism $\pi :C\to\Pone$ of degree $d$ such that $\pi
^*(0)=d(x)$, $\pi ^*(\infty )=d(y)$ and $\pi$ does not factor through a cover
$z\in\Pone\mapsto z^r\in\Pone$ for some $r>1$. The previous lemma shows that
all such covers are of the same topological type. This
implies that $Y_d$ is irreducible. So every irreducible component of $Y$ is
equal to some $Y_d$.

Let $\JJ_g\to\M _g$ be the universal Jacobian and let $q:\CC _g^2\to\JJ
_g$ be the Abel-Jacobi map $(C,x,y)\mapsto (x)-(y)\in J(C)$. Then
$Y_d=q^{-1}\JJ
_g\la d\ra $. Since $Y_d$ has the correct codimension $g$ in $\CC _g^2$, it
follows that $[Y_d]$ is a positive multiple of $q^*[\JJ _g\la d\ra ]$.
According
to \refer{2.9}, $[\JJ _g\la d\ra ]$ is a positive multiple of the class of the
zero section in $A^g(\JJ _g )$ and so the proof is complete.
\enddemo

\demo{Proof of \refer{1.2}} First observe that the direct image of
$K_1^{1+d_1}\cdots K_n^{1+d_n}$ under the forgetful morphism $\CC _g^n\to\M _g$
equals $\kappa _{d_1}\cdots\kappa _{d_n}$. Now the direct image of the class of
an irreducible component of $X^k$ of codimension $k$ under $\CC _g^n\to\M _g$
is
zero unless the image has the correct codimension $k-n$. In particular, a
nonzero image requires $k\ge n$. It follows that any product in the
tautological
classes of degree $d$ can be represented by a linear combination of the
irreducible components of the locus in $\M _g$ that parametrizes the curves $C$
that admit a covering $\pi :C\to \Pone$ of degree $\le 2g-2$ totally ramified
over $\infty$ and with at most $g-1-d$ points over $0$. The rest follows
immediately from \refer{1.1}.
\enddemo

\enddocument